\documentclass[11pt]{article}
\input{epsf.sty}
\usepackage{hyperref}
\usepackage{amsmath, amssymb}
\newtheorem {thm}{Theorem}[section]
\newtheorem {prop}[thm]{Proposition}
\newtheorem {lem}[thm]{Lemma}
\newtheorem {cor}[thm]{Corollary}
\newtheorem {defn}[thm]{Definition}

\def\Cox{\hfill \Box}
\def\Z{{\Bbb Z}}
\def\R{{\Bbb R}}
\def\P{{\Bbb P}}

\def\0{{\bf 0}}

\def\eps{\varepsilon}
\def\dist{{\rm dist}}
\def\parent{{\rm parent}}

\def\b{\beta}

\def\d{\delta}

\def\phi{\varphi}

\def\l{\lambda}

\def\L{\Lambda}

\def\T{\T}

\def\PP{{\cal P}}
\def\Norm{Norm.\,}

\begin{document}

\title{Gibbs properties of the fuzzy Potts model \\
on trees and in mean field}
\author{Olle H\"aggstr\"om\thanks{Research supported by the 
Swedish Research Council} \footnote{Mathematical Statistics, 
Chalmers University of Technology,
412 96 G\"oteborg,
SWEDEN,
\texttt{olleh@math.chalmers.se},
\texttt{http://www.math.chalmers.se/\~{ }olleh/}} \, and Christof 
K\"ulske\footnote{Weierstrass-Institut f\"ur 
Angewandte Analysis und Stochastik,
Mohrenstr. 39,
D-10117 Berlin,
GERMANY,
\texttt{kuelske@wias-berlin.de},
\texttt{http://www.wias-berlin.de/people/kuelske/} }}
\maketitle

\begin{abstract}
We study Gibbs properties of the fuzzy Potts model in the mean field
case (i.e.\ on a complete graph) and on trees. For the mean field case,
a complete characterization of the set of temperatures for which 
non-Gibbsianness happens is given. 
The results for trees are somewhat less explicit, but we do show for general
trees that non-Gibbsianness of the fuzzy Potts model happens exactly
for those temperatures where the underlying Potts model has multiple
Gibbs measures. 
\end{abstract}

\section{Introduction} \label{sect:intro}

It used to be taken for granted that simple transformations of Gibbs measures
are themselves Gibbsian. A few counterexamples were found in the 70's
and 80's \cite{GP,S}, but these were usually referred to as being
somehow exceptional or pathological. In the seminal paper from 1993
by van Enter, Fern\'andez and Sokal \cite{vEFS}, further examples were
found, and a systematic study of Gibbsianness vs.\ non-Gibbsianness
of large classes of transformed or projected versions of Gibbs systems began;
see \cite{vEFK,MVV,vE,FP,K1,vEMS,vEMSS,K2,vEFdHR,H3} for some of the subsequent
work in this area. 

In particular, Gibbs properties of the so-called fuzzy Potts model
were studied in Maes and Vande Velde \cite{MVV} and H\"aggstr\"om \cite{H3}.
Like almost all work in the study of Gibbsianness vs.\ non-Gibbsianness,
these papers focused on the case where the underlying lattice is $\Z^d$.
Two exceptions are H\"aggstr\"om \cite{H96} and K\"ulske \cite{K3} 
where these issues are
studied for certain models living on trees and on complete graphs
(known as the Curie--Weiss or mean field case), 
respectively. In this paper, we shall
continue in the directions of \cite{H96} and \cite{K3} by 
studying Gibbs properties of the fuzzy Potts model on trees and in
the mean field setup. 

The fuzzy Potts model arises, loosely speaking, from the standard 
$q$-state Potts 
model by looking at it with a pair of glasses that prevents from 
distinguishing some of the spin values; see Section \ref{sect:model}
for precise definitions. 
This makes the fuzzy Potts model one of the most basic examples of a
hidden Markov random field \cite{KGK}, and it has also turned out
to be useful in the study of percolation-theoretic properties
of the underlying Potts model \cite{C,H1}. 
Maes and Vande Velde \cite{MVV} speculated
that Gibbsianness of the fuzzy Potts model on $\Z^d$
might hold precisely in the Gibbs uniqueness regime (i.e., above the
critical temperature) of the underlying Potts model, but
this was shown in \cite{H3} not to be the case: non-Gibbsianness
of the fuzzy Potts model happens also for some parameter values
where the underlying Potts model has a unique Gibbs measure. In the
following result, which is our main result for trees, we see that
the desired equivalence between on one hand
Gibbsianness of the fuzzy Potts model and on the other hand
Gibbs uniqueness of the underlying Potts model does hold when
$\Z^d$ is replaced by a tree. 
\begin{thm}  \label{thm:main_result_on_trees}
Consider the $q$-state Potts model on a tree $\Gamma$ at inverse 
temperature $\beta$, and let $s$ and $r_1, \ldots, r_s$ be positive
integers with $1<s<q$ and $\sum_{i=1}^s r_i =q$. The set ${\cal G}$
of Gibbs measures for this Potts model contains an element
whose corresponding fuzzy Potts
measure with spin partition $(r_1, \ldots, r_s)$ is non-Gibbsian,
if and only if
$| {\cal G} | > 1$. 
\end{thm}
Here $| \cdot |$ denotes cardinality, while the remaining notation and
terminology will be explained in later sections. 

We move on to the mean-field fuzzy Potts model, which lives on
a complete graph on $N$ vertices and for which we consider asymptotics as 
$N\rightarrow \infty$. Here
the situation is quite different. 
Before we state our main result a few general remarks are in order. 
First of all one has to be careful to find the right way of asking for  
``Gibbssianness'' vs. ``non-Gibbsianness'' for mean-field models. It must 
be asked in  an appropriate sense
if we want  to see non-trivial behavior that reflects the 
lattice-phenomenon in a natural way. 
We remind the reader that the Gibbs measures of simple
mean-field models usually converge
weakly to (linear combinations of) product measures.
A (non-trivial) linear
combination of product measures is non-Gibbsian and has each spin configuration
as a point of discontinuity when we are looking at it in the {\it product topology}  \cite{EL96}. 
So the problem of finding non-Gibbsianness in mean-field models would always 
have a trivial (negative) answer as soon as there is a phase transition, 
and a trivial (positive) answer as soon as there is no phase transition, 
independently of the model.  
We stress that this is a very different phenomenon than the one happening for the 
fuzzy Potts model on the tree described above. 
However, if we don't want to stop at this point but want to see something meaningful 
we must proceed differently. 
As it was argued in \cite{K3} {\it non-Gibbsianness for mean-field models} should be understood as {\it 
discontinuity
of conditional probabilities as a function of the conditioning}, 
but the notion of continuity must {\it not} be taken with respect to
product topology.  More precisely, 
we need to perform the following limiting procedure. 
\begin{enumerate}
\item
Take the conditioning of the conditional probabilities 
of the finite volume Gibbs-measures while staying in
finite volume. Due to permutation invariance, 
these conditional probabilities are automatically 
volume-dependent {\it functions of the empirical average}
over all the spins in the conditioning. 
\item
Show that  the large volume-limit 
for these functions exists, and look at their continuity properties. 
\end{enumerate}

When these limiting conditional probabilities 
 are discontinuous, we have found an analogue of ``non-Gibbsian'' behavior 
in the mean-field model.  
When they are continuous, the mean-field model behaves in a ``Gibbsian'' way. 
In the case of non-Gibbsian behavior 
we can carry the analogue between mean-field and lattice 
to the notion of ``almost sure Gibbsianness'' (that is familiar on 
the lattice). For the mean-field model we look  
at the size of the set of the discontinuity points in the large volume-limit,  with respect to the 
limiting measure on the empirical distribution. If the discontinuity points get measure 
zero, we have found the mean-field analogue of ``almost sure Gibbsian'' behavior. 

An analysis of this sort was carried out in \cite{K3} for the decimation transformation of the 
Ising ferromagnet, 
and examples of joint measures in random systems including 
the random field Ising model. For the models we were looking at  we saw a surprising analogy 
between mean-field and lattice results. 

We are now ready to state our main result for the mean-field version of the fuzzy Potts 
model in short form. Precise definitions and more details will be given 
in  \ref{sect:mean-field}.

\begin{thm}  \label{thm:main_result_on_complete_graph}
Consider the $q$-state mean-field Potts model  at inverse 
temperature $\beta$, and let $s$ and $r_1, \ldots, r_s$ be positive
integers with $1<s<q$ and $\sum_{i=1}^s r_i =q$. 
Consider the limiting conditional probabilities of the 
corresponding fuzzy Potts
model with spin partition $(r_1, \ldots, r_s)$. 
\begin{description}
\item {\bf (i)} Suppose that $r_i\leq 2$ for all $i=1,\dots,s$. 
Then the limiting conditional probabilities are continuous functions 
of the empirical mean of the conditioning, for all $\b\geq 0$. 
\end{description}
Assume that $r_i\geq 3$ for some $i$ and put
$r_*:=\min\{ r\geq 3, r=r_i \hbox{ for some }i=1,\dots,s\}$.  
Denote by $\b_c(r)$ the inverse critical temperature of 
the $r$-state Potts model. Then the following holds. 
\begin{description}
\item{\bf (ii)} 
The limiting conditional probabilities are continuous for all  
$\b<\b_c(r_*)$.
\item{\bf (iii)} 
The limiting conditional probabilities are discontinuous for all  
$\b\geq \b_c(r_*)$.  
\item{\bf (iv)} The set of discontinuity points has zero measure 
in the infinite volume limit in all cases. 
\end{description}  
\end{thm}
Thus, we have a rather complete picture for the limiting behavior of the model 
on complete graphs. Note that from (iii) follows in particular 
that there is an interesting range of temperatures $\b_c(r_*)\leq \b<\b_c(q)$ 
when the underlying Potts model shows no phase transition 
but the fuzzy model is non-Gibbsian. (It is well-known that $\b_c(q)$ 
is increasing with $q$.) 
As mentioned above, the existence of such a region was shown 
on the lattice in \cite{H3}; in the present mean-field model 
the lower endpoint of the interval is moreover proved to be $\b_c(r_*)$ 
(which is only a conjecture on the lattice). 
For such a non-Gibbsianness to occur in mean-field we need  however 
that there is at least one fuzzy class containing 
three or more spin-values. This is due to the fact that the discontinuity 
of the limiting conditional 
probabilities is related to a first order transition within 
one fuzzy class, and such a 
transition exists if and only if there are at least three spin values. 

Controlling the size of the set of discontinuities 
is a more subtle task, but we manage in Theorem 
\ref{thm:main_result_on_complete_graph} (iv) to provide the complete answer
in the mean-field case: almost sure Gibbsianness holds regardless of the
choice of parameter values. 

The rest of the paper is organized as follows. 
In Section \ref{sect:model} we define the models. 
In Section \ref{sect:quasilocal} we briefly explain why, in the case of the
fuzzy Potts model, Gibbsianness is the same thing as so-called quasilocality.
Our main results for trees are stated and proved in Section 
\ref{sect:trees}, whereas those in the mean field setup are
treated in Section \ref{sect:mean-field}.  We mention that Section 
\ref{sect:mean-field} can be read independently of Sections
\ref{sect:quasilocal} and \ref{sect:trees}.

\section{The models}  \label{sect:model}

In this section we give the definitions (following \cite{H3})
of the Potts model and the fuzzy
Potts model, first on finite graphs, and then on infinite graphs.
The results in Section \ref{subsect:Potts_infinite}
concerning infinite-volume limits of the Potts model
date back to Aizenman et al.\ \cite{ACCN}; see also \cite{GHM} for a
detailed account of these results. 

\subsection{Potts in finite volume}

For a positive integer $q$,  
the $q$-state {\bf Potts model} on a finite graph $G=(V,E)$ is a random
assignment of $\{1, \ldots, q\}$-valued spins to the vertices of $G$. The
Gibbs measure $\pi^G_{q,\beta}$ for the $q$-state Potts model on $G$ 
at inverse temperature $\beta\geq 0$, is
the probability measure $\pi^G_{q,\beta}$ on $\{1, \ldots, q\}^V$ 
which to each $\xi\in \{1, \ldots, q\}^V$ assigns probability
\begin{equation} \label{eq:Potts}
\pi^G_{q,\beta}(\xi)=\frac{1}{Z^G_{q, \beta}}\exp\left(
2\beta \sum_{\langle x,y \rangle\in E} I_{\{\xi(x) = \xi(y)\}} \right) \, .
\end{equation}
Here $\langle x,y \rangle$ denotes the edge connecting $x,y \in V$, $I_A$ 
is the indicator function of the event $A$, and $Z^G_{q, \beta}$ is a
normalizing constant.

\subsection{Fuzzy Potts in finite volume}

Next, let $s$ and $r_1, \ldots, r_s$ be positive 
integers such that $\sum_{i=1}^s r_i = q$. The {\bf fuzzy Potts model}
on $G$ with these parameters 
arises by taking the $q$-state Potts model on $G$, and then identifying
the first $r_1$ Potts states with a single fuzzy spin value $1$, 
the next $r_2$ of the states
with fuzzy spin value $2$, and so on. A more precise definition is as follows. 
Fix $q$, $\beta$ and $(r_1, \ldots, r_s)$ as above. Let $X$ be a
$\{1, \ldots, q\}^V$-valued random object distributed according to the Gibbs
measure $\pi^G_{q, \beta}$. Then take $Y$ to be the 
$\{1\ldots, s\}^V$-valued random object obtained from $X$ by setting
\begin{equation} \label{eq:def_of_fuzzy_Potts}
Y(x) = \left\{ 
\begin{array}{cl}
1 & \mbox{if } X(x) \in \{1, \ldots, r_1\} \\
2 & \mbox{if } X(x) \in \{r_1 +1, \ldots, r_1 + r_2\} \\
\vdots & \hspace{7mm} \vdots \\
s & \mbox{if } X(x) \in \{ q-r_s + 1, \ldots, q\} 
\end{array} \right. 
\end{equation}
for each $x \in V$. 
We write $\mu^G_{q, \beta, (r_1, \ldots, r_s)}$ for the probability
measure on $\{1, \ldots, s\}^V$ which describes the distribution of $Y$, 
and call it the fuzzy Potts measure with
parameters $q$, $\beta$, and $(r_1, \ldots, r_s)$. We call
$(r_1, \ldots, r_s)$ the {\bf spin partition} for this fuzzy Potts model.

Of course, $\mu^G_{q, \beta, (r_1, \ldots, r_s)}$ is uninteresting
for $s=1$, whereas for $s=q$ it just reproduces the ordinary Potts model. We 
therefore require that $1<s<q$, and consequently that $q\geq 3$. 

\subsection{Potts in infinite volume} \label{subsect:Potts_infinite}

Now let $G=(V,E)$ be infinite and locally finite. For
$W \subset V$, we define its boundary $\partial W$ as
\[
\partial W = \{x \in V \setminus W: \exists y \in W \mbox{ such that }
\langle x,y \rangle \in E\} \, .
\]
A probability measure $\pi$ on $\{1, \ldots, q\}^V$ is said to be
a Gibbs measure for the $q$-state Potts model on $G$ at inverse temperature
$\beta$, if it admits conditional probabilities such that for all finite
$W \subset V$, all $\xi \in \{1, \ldots, q\}^W$ and all
$\eta \in \{1, \ldots, q\}^{V \setminus W}$ we have
\begin{eqnarray} \nonumber 
\lefteqn{\pi(X(W)=\xi \, | \, X(V \setminus W)= \eta)} \hspace{10 mm} \\ 
& = & \frac{1}{Z^{W, \eta}_{q, \beta}}
\exp \left( 2 \beta \left( \sum_{\langle x,y\rangle \in E \atop x,y \in W}
I_{\{\xi(x) = \xi(y)\}} + 
\sum_{\langle x,y\rangle \in E \atop x\in W, y \in \partial W}
I_{\{\xi(x) = \eta(y)\}} \right) \right)
\label{eq:def_of_Potts}
\end{eqnarray}
where the normalizing constant $Z^{W, \eta}_{q, \beta}$ 
depends on $\eta$ but not on $\xi$. Note that the corresponding
relation holds in the finite graph case where $\pi$ is defined by 
(\ref{eq:Potts}). 

The basic examples of Gibbs measures for the Potts model are constructed as
follows. Let $\Lambda= \{\Lambda_n\}_{n=1}^\infty$ denote a sequence 
of subsets of $V$, which is an exhaustion of
$V$ in the sense that
(i) each $\Lambda_n$ is finite,
(ii) $\Lambda_1 \subset \Lambda_2 \subset \cdots$, and
(iii) $\bigcup_{n=1}^\infty \Lambda_n = V$.
Let $G_n$ denote
the graph whose vertex set is $\Lambda_n \cup \partial \Lambda_n$, and whose
edge set consists of pairs of vertices in $\Lambda_n \cup \partial \Lambda_n$
at 
distance $1$ from each other. It is well-known
that the Gibbs measures $\pi^{G_n}_{q, \beta}$ converge to a probability
measure on $\{1, \ldots, q\}^V$ which is a Gibbs measure for the Potts
model on $G$ with the given parameters. Convergence takes place in the
sense that probabilities of cylinder sets converge. The limiting probability
measure on $\{1, \ldots ,q\}^V$ is
denoted $\pi^{G, 0}_{q, \beta}$, and is called the Gibbs measure
(for the Potts model on $G$ with the given parameters) with
{\bf free boundary condition}. Other Gibbs measures are those with
so-called {\bf spin $i$ boundary condition}, denoted 
$\pi^{G, i}_{q, \beta}$, for $i=1, \ldots q$. These are obtained by
conditioning $\pi_{q, \beta}^{G_n}$ on taking spin value $i$ all over
$\partial \Lambda_n$
and then taking limits as $n \rightarrow\infty$. 
The existence of these limits, and the fact that
each of the measures 
$\pi^{G, 0}_{q, \beta}, \ldots, \pi^{G, q}_{q, \beta}$
is independent of the particular choice of 
exhaustion $\{\Lambda_n\}_{n=1}^\infty$,
follows from the work of Aizenman et al.\ \cite{ACCN}.

The Gibbs measures 
$\pi^{G, 0}_{q, \beta}, \pi^{G, 1}_{q, \beta}, \ldots
\pi^{G, q}_{q, \beta}$
may or may not coincide depending on $G$ and on
the parameter values. It is a fundamental result from 
\cite{ACCN} that the occurence of more than one distinct
Gibbs measure is (for fixed $G$ and $q$) increasing in $\beta$. Hence,
there exists a critical value
$\beta_c=\beta_c(G,q)\in (0,\infty)$, such that for $\beta< \beta_c$, there 
is only one Gibbs measure (so that in particular
$\pi^{G, 0}_{q, \beta}= \cdots = \pi^{G, q}_{q, \beta}$), whereas
for $\beta>\beta_c$, there are multiple Gibbs measures (and moreover the
measures $\pi^{G, 0}_{q, \beta}, \ldots, \pi^{G, q}_{q, \beta}$ are
all different). The critical value may be $\infty$ if the graph is 
``too small'' or $0$ if the graph is ``too large'' (requiring
unbounded degree and more than that) but in many interesting cases
there is a nontrivial critical value $\beta_c \in (0, \infty)$, such as for
cubic lattices in $d\geq 2$ dimensions and regular trees of degree
at least $3$. Yet another important result from \cite{ACCN} is that
nonuniqueness of Gibbs measures is equivalent to having
\begin{equation}  \label{eq:positive_magnetization}
\pi^{G, 1}_{q, \beta} (\mbox{spin $1$ at $x$}) > \frac{1}{q}
\end{equation}
for some $x\in V$, and that if $G$ is connected, then this is in turn
equivalent to having (\ref{eq:positive_magnetization}) for {\em every}
$x\in V$. (For symmetry reasons, we have
\begin{equation}  \label{eq:zero_magnetization}
\pi^{G, 0}_{q, \beta} (\mbox{spin $1$ at $x$}) = \frac{1}{q}
\end{equation}
for every $x\in V$. Whenever we are in the uniqueness regime of the parameter
space, we then of course have (\ref{eq:zero_magnetization}) with
$\pi^{G, 0}_{q, \beta}$ replaced by any of the other Gibbs measures
$\pi^{G, i}_{q, \beta}$.)

\subsection{Fuzzy Potts in infinite volume}

Given the Gibbs measures 
$\pi^{G, 0}_{q, \beta}, \pi^{G, 1}_{q, \beta}, \ldots,
\pi^{G, q}_{q, \beta}$, we define
fuzzy Potts measures as in the case of finite graphs. More precisely, 
for $q$, $\beta$, and $(r_1, \ldots, r_s)$ as above, and
$i\in \{0, \ldots, q\}$, we define the fuzzy Potts measure 
$\mu^{G,i}_{q, \beta, (r_1, \ldots, r_s)}$ to be the distribution of
the $\{1, \ldots, s\}^V$-valued random object $Y$ obtained by first
picking $X\in \{1,\ldots, q\}^V$ according to the Gibbs measure
$\pi^{G, i}_{q, \beta}$, and then constructing $Y$ from $X$ as in
(\ref{eq:def_of_fuzzy_Potts}). 

\section{Gibbsianness and quasilocality}  \label{sect:quasilocal}

When $S$ is a finite set, $G=(V,E)$ is an infinite locally finite graph, and
$\mu$ is a probability measure on $S^V$, it is well known 
(see, e.g., \cite[Thm.\ 2.12]{vEFS}) that
$\mu$ is Gibbsian if and only if it satisfies the properties of
quasilocality and uniform nonnullness. The latter property
means that $\mu$ admits conditional probabilities such that
\[
\min_{s\in S} \inf_{\eta \in S^{V \setminus \{x\}}}
\mu(X(x)=s \, | \, X(V \setminus \{x\}) = \eta) > 0
\]
for each $x\in V$. Uniformly nonullness holds in the Potts model, and it
is easy to see that this property is inherited by the fuzzy Potts model;
see \cite[Lem.\ 4.5]{H3}. Hence, the problem of determining whether
the fuzzy Potts model with given parameter values is Gibbsian is reduced
to that of whether it is quasilocal. Quasilocality is defined as follows,
where (as in Section \ref{sect:model})
$\Lambda= \{\Lambda_n\}_{n=1}^\infty$ is an exhaustion of $V$
(the definition does not depend on the particular choice of $\Lambda$). 
\begin{defn}  \label{defn:quasilocality}
Let $S$ be a finite set and let
$G=(V,E)$ be an infinite locally finite graph. A probability measure $\mu$
on $S^V$ is said to be {\bf quasilocal} if it admits conditional 
probabilities such that for all finite $W \subset V$ and all 
$\xi \in S^W$ we have
\begin{equation}  \label{eq:def_quasilocality}
\lim_{n \rightarrow\infty} 
\sup_{\eta, \eta' \in S^{V \setminus W} \atop \eta(\Lambda_n \setminus W)
=\eta'(\Lambda_n \setminus W)} 
\Big| \mu(X(W)=\xi \, | \, X(V \setminus W)=\eta) - 
\mu(X(W)=\xi \, | \, X(V \setminus W)=\eta') \Big| = 0 \, .
\end{equation}
\end{defn}
Because of the asserted equivalence between Gibbsianness and quasilocality
for the fuzzy Potts model, we shall in the following focus
entirely on quasilocality. In Section \ref{sect:trees} on trees, this
means studying the property in Definition \ref{defn:quasilocality}
verbatim, whereas in Section \ref{sect:mean-field} we need to adapt
the definition of quasilocality somewhat (following \cite{K3}), as
hinted in Section \ref{sect:intro}.

\section{The fuzzy Potts model on trees}  \label{sect:trees} 

\subsection{Trees: definitions}

A tree $\Gamma$ is a connected graph without cycles. In addition to
these properties, we assume that $\Gamma$ is locally finite, and we
denote its vertex set and edge set by $V_\Gamma$ and $E_\Gamma$, respectively. 
Pick an arbitrary vertex in $\rho\in V_\Gamma$ and call it the 
{\bf root} of $\Gamma$. For $x,y\in V_\Gamma$, let $\dist(x,y)$ denote
the graph-theoretic distance between $x$ and $y$ in $\Gamma$. 
If $x$ and $y$ share an edge and $\dist(y,\rho)=\dist(x,\rho)+1$, then
we call $y$ a {\bf child} of $x$, and $x$ is the {\bf parent} of $y$. 
More generally, if $x$ is on the unique self-avoiding path from $\rho$
to $y$, then $y$ is called a {\bf descendant} of $x$, and $x$ is an
{\bf ancestor} of $y$. Each vertex $x$
except for the root has exactly one parent, denoted $\parent(x)$
while the number of children may vary. If two vertices $x$ and $y$ 
have the same parent, then we call them {\bf siblings}. 

An important example is when, for some $d \geq 2$, the root has 
$d+1$ children and all others have $d$ children; this is referred to as the
regular tree with degree $d$. See, e.g., \cite{P} for a variety of
other interesting examples of trees. 

For $n=0,1,\ldots$, let $\Gamma_n = (V_{\Gamma_n}, E_{\Gamma_n})$ be
the subgraph (subtree) of $\Gamma$ given by
\[
V_{\Gamma_n} = \{ x \in V_\Gamma : \, \dist(x,\rho) \leq n\} \, 
\]
and
\[
E_{\Gamma_n} = \{ e \in E_{\Gamma} : \mbox{both endpoints of $e$ are
in }V_{\Gamma_n}\} \, ,
\]
and note that $\{V_{\Gamma_n}\}_{n=1}^\infty$ is an exhaustion of $V_\Gamma$.
For $x\in \Gamma$, let $\Gamma_{(x)}$ denote the induced subtree of $\Gamma$
whose vertex set consists of $x$ and all its descendants. In other words,
$\Gamma_{(x)}= (V_{\Gamma_{(x)}}, E_{\Gamma_{(x)}})$ with
\[
V_{\Gamma_{(x)}} = \{y \in V_{\Gamma} : \, x \mbox{ is an ancestor of } y \} 
\]
and
\[
E_{\Gamma_{(x)}} = \{e \in E_{\Gamma} : \mbox{both endpoints of $e$ are
in } V_{\Gamma_{(x)}} \} \, .
\]
Finally, for $x\in \Gamma$ and $n \geq \dist(x,\rho)$, define the subtree
$\Gamma_{(x,n)}= (V_{\Gamma_{(x,n)}}, E_{\Gamma_{(x,n)}})$ by setting
\[
V_{\Gamma_{(x,n)}} = V_{\Gamma_{(x)}} \cap V_{\Gamma_n}
\]
and
\[
E_{\Gamma_{(x,n)}} = E_{\Gamma_{(x)}} \cap E_{\Gamma_n} \, .
\]

\subsection{Proofs}

The key results for proving Theorem \ref{thm:main_result_on_trees} are
the following two propositions.
\begin{prop}  \label{prop:free}
Let $\Gamma$ be a tree, and fix the parameter values $q$,
$\beta$, $s$ and $(r_1, \ldots r_s)$ with $1<s<q$ for the Potts model
and the fuzzy Potts model on $\Gamma$. Then the fuzzy Potts measure
$\mu^{G,0}_{q, \beta, (r_1, \ldots, r_s)}$ corresponding to the
Gibbs measure with free boundary condition, is quasilocal. 
\end{prop}
\begin{prop}  \label{prop:wired}
Let $\Gamma$ be a tree, and fix the parameter values $q$,
$\beta$, $s$ and $(r_1, \ldots r_s)$ with $1<s<q$ and $r_1 > 1$
for the Potts model
and the fuzzy Potts model on $\Gamma$. Suppose that
$\pi^{G, 1}_{q, \beta} \neq \pi^{G, 0}_{q, \beta}$. Then
$\mu^{G,1}_{q, \beta, (r_1, \ldots, r_s)}$ is nonquasilocal. 
\end{prop}
{\bf Proof of Theorem \ref{thm:main_result_on_trees} from Propositions
\ref{prop:free} and \ref{prop:wired}:}
Since $s<q$, we must have $r_i>1$ for some $i\in \{1, \ldots, s\}$,
and there is no loss of generality in assuming that $r_1>1$. If
$| {\cal G}|>1$, then
$\pi^{G,1}_{q, \beta} \neq \pi^{G,0}_{q, \beta}$ due to 
(\ref{eq:positive_magnetization}) and (\ref{eq:zero_magnetization}).
Hence, using Proposition \ref{prop:wired}, 
$\mu^{G,1}_{q, \beta, (r_1, \ldots, r_s)}$ is nonquasilocal
and therefore non-Gibbsian, and the `if' part of the
theorem is established. For the `only if' part, note that
if $| {\cal G}|=1$, then ${\cal G}=\{\pi^{G,0}_{q, \beta}\}$, so
that $\mu^{G,0}_{q, \beta, (r_1, \ldots, r_s)}$ is the only fuzzy Potts
measure, which by Proposition \ref{prop:free} is quasilocal
and therefore Gibbsian. $\Cox$

\medskip\noindent
It remains to prove Propositions \ref{prop:free} and \ref{prop:wired}.
To this end, we need to introduce the notion of a tree-indexed Markov chain
on $\Gamma$, and its relation to 
Gibbs measures for the Potts model on $\Gamma$. 
This relation is well-known for regular trees (see for instance 
\cite{Sp,Z,BW}), while the extension to general trees seems
to be less well-studied. 

Let $(x_0, x_1, \ldots)$ be an enumeration of $V_{\Gamma}$ such that
the root $\rho$ comes first $(x_0=\rho$), then all
vertices in $V_{\Gamma_1} \setminus \{\rho\}$, then all vertices in 
$V_{\Gamma_2} \setminus V_{\Gamma_1}$, and so on. 
Fix $q$, let $\nu$ be a probability measure on $\{1, \ldots, q\}$
(which will play the role of an initial distribution), and let
$P = (P_{ij})_{i,j \in \{1, \ldots, q\}}$ be a transition matrix.
Let $X$ be the $\{1, \ldots, q\}^{V_{\Gamma}}$-valued random spin
configuration obtained as follows. First pick 
$X(x_0)\in \{1, \ldots, q\}$
according to $\nu$. Then, inductively, once $X(x_0), \ldots, X(x_n)$ have 
been determined, pick $X(x_{n+1})\in \{1,\ldots, q\}$ with distribution
$(P_{i1}, \ldots, P_{iq})$ where $i=X(\parent(x_{n+1}))$. For obvious
reasons, $X$ is called a tree-indexed Markov chain on $\Gamma$. 

There is sometimes reason to consider inhomogeneous tree-indexed
Markov chains, where the transition matrix $P$ is allowed to depend
on where in the tree we are: for every $x\in V_\Gamma \setminus \{ \rho \}$,
we then 
have a transition matrix $P^x = (P^x_{ij})_{i,j \in \{1, \ldots, q\}}$,
and $X$ is generated as above with $X(x)$ chosen according to 
the distribution $(P^x_{i1}, \ldots, P^x_{iq})$ where 
$i=X(\parent(x))$. 

It is readily checked that a (possibly inhomogeneous)
tree-indexed Markov chain $X$ is also a Markov
random field on $\Gamma$, meaning that for any finite $W \subset V_\Gamma$,
the conditional distribution of $X(W)$ given $X(V_\Gamma \setminus W)$ depends
on $X(V_\Gamma \setminus W)$ only via  $X(\partial W)$. Hence the supremum in
(\ref{eq:def_quasilocality}) becomes $0$ for all $n$ large enough so that
$\Lambda_n$ contains $W \cap \partial W$, so that we have the following
lemma.
\begin{lem}  \label{lem:chains_are_quasilocal}
The distribution of any homogeneous or inhomogeneous
tree-indexed Mar\-kov chain on $\Gamma$ is quasilocal.
\end{lem}
Fix $\beta\geq 0$, and consider the tree-indexed Markov chain
given by $\nu=(\frac{1}{q}, \ldots, \frac{1}{q})$ and transition matrix
$P=(P_{ij})_{i,j \in \{1, \ldots, q\}}$ given by
\begin{equation}  \label{eq:transition_matrix}
P_{ij} = \left\{
\begin{array}{ll}
\frac{e^{2\beta}}{e^{2\beta}+q-1} & \mbox{if } i =j \\
\frac{1}{e^{2\beta}+q-1} & \mbox{otherwise.}
\end{array} \right. 
\end{equation}
Let $X\in\{1, \ldots, q\}^{V_{\Gamma}}$ be given by this particular
tree-indexed Markov chain. By directly checking (\ref{eq:Potts}),
we see that $X(\Lambda_n)$ has distribution $\pi^{\Gamma_n}_{q,\beta}$.
By taking limits as $n\rightarrow \infty$ and considering the
construction of $\pi^{G,0}_{q,\beta}$ in Section \ref{subsect:Potts_infinite},
we see that $X$ is distributed according to the Gibbs measure
$\pi^{\Gamma,0}_{q,\beta}$ for the Potts model on $\Gamma$ with free boundary
condition. 

\medskip\noindent
{\bf Proof of Proposition \ref{prop:free}:}
Construct $X\in \{1, \ldots, q\}^{V_{\Gamma}}$ sequentially as above,
with $\nu= (\frac{1}{q}, \ldots, \frac{1}{q})$ and $P$ given by
(\ref{eq:transition_matrix}), and let
$Y\in \{1, \ldots, r\}^{V_{\Gamma}}$ from $X$ as in 
(\ref{eq:def_of_fuzzy_Potts}). Then the conditional distribution
of $Y(x_{n+1})$ given $X(x_0), \ldots, X(x_n)$ such that 
$X(\parent(x_{n+1}))=i$ and $Y(\parent(x_{n+1}))=k)$, is given by
\begin{equation}  \label{eq:Y_is_a_Markov_chain}
\P(Y(x_{n+1}) = l \, | \, \cdots) = \left\{
\begin{array}{ll}
\frac{e^{2\beta}+r_k-1}{e^{2\beta}+q-1} & \mbox{if } l=k \\
\frac{r_k}{e^{2\beta}+q-1} & \mbox{otherwise,}
\end{array} \right. 
\end{equation}
which follows by summing over the possible values of $X(x_{n+1})$. 
Note that the right-hand side of (\ref{eq:Y_is_a_Markov_chain}) depends on
$X(x_0), \ldots, X(x_n)$ only through $Y(\parent(x_{n+1}))$. It follows
that $Y$ is a tree-indexed Markov chain with state space
$\{1, \ldots, s\}$, initial distribution $(\frac{r_1}{q}, \ldots, 
\frac{r_s}{q})$ and transition matrix $P=(P_{kl})_{k,l \in \{1, \ldots, s\}}$
given by
\begin{equation}  \label{Markov_chain_for_Y_with_free_bc}
P_{kl} = \left\{
\begin{array}{ll}
\frac{e^{2\beta}+r_l-1}{e^{2\beta}+q-1} & \mbox{if } l=k \\
\frac{r_l}{e^{2\beta}+q-1} & \mbox{otherwise.}
\end{array} \right. 
\end{equation} 
Quasilocality of $Y$ now follows from Lemma \ref{lem:chains_are_quasilocal}.
$\Cox$

\medskip\noindent
For the proof of Proposition \ref{prop:wired}, we need to consider the
tree-indexed Markov chain on $\Gamma$ corresponding to the Gibbs measure
$\pi^{\Gamma,1}_{q, \beta}$ with the ``all $1$'' boundary condition. This 
is a bit more complicated than the case of $\pi^{\Gamma,0}_{q, \beta}$ due
to the lack of full symmetry among the spin values. 

For $x\in V_{\Gamma}$, consider the Gibbs measure 
$\pi^{\Gamma_{(x)},1}_{q, \beta}$, and in particular the probability
$\pi^{\Gamma_{(x)},1}_{q, \beta}(\mbox{spin $1$ at } x)$, which we denote
by $a_x$. (Note that $a_x$ is in general distinct from 
$\pi^{\Gamma,1}_{q, \beta}(\mbox{spin $1$ at } x)$, because it fails to
take into account, e.g., the possible influence from 
$\parent(x)$ on $x$.) For symmetry reasons, the 
$\pi^{\Gamma_{(x)},1}_{q, \beta}$-distribution of
the spin at $x$ is 
\[
\left(a_x, \frac{1-a_x}{q-1}, \frac{1-a_x}{q-1}, \ldots, \frac{1-a_x}{q-1}
\right) \, .
\]
Also define
\begin{equation} \label{eq:def_b_x}
b_x = \frac{a_x}{(1-a_x)/(q-1)} = 
\frac{\pi^{\Gamma_{(x)},1}_{q, \beta}(\mbox{spin $1$ at } x)}{
\pi^{\Gamma_{(x)},1}_{q, \beta}(\mbox{spin $2$ at } x)} \, .
\end{equation}
The constants $\{b_x\}_{x\in V_\Gamma}$ satisfy the following recursion. 
\begin{lem}  \label{lem:wired_recursion}
Suppose $x \in V_\Gamma$ is a vertex with $k$ children $y_1 , \ldots, y_k$. 
We then have
\begin{equation} \label{eq:wired_recursion}
b_x = \frac{\prod_{i=1}^k (e^{2\beta} b_{y_i} + q - 1)}{
\prod_{i=1}^k (e^{2\beta} + b_{y_i} + q - 2)} \, .
\end{equation} 
\end{lem}
{\bf Proof:}
For $n$ large enough so that $x \in V_{\Lambda_n}$, 
define, as a finite-volume analogue of (\ref{eq:def_b_x}), 
\[
b_{x,n} = \frac{\pi^{\Gamma_{(x,n)},1}_{q, \beta}(\mbox{spin $1$ at } x)}{
\pi^{\Gamma_{(x,n)},1}_{q, \beta}(\mbox{spin $2$ at } x)} \, ,
\]
where $\pi^{\Gamma_{(x,n)},1}_{q, \beta}$ is the finite-volume Gibbs measure
for $\Gamma_{(x,n)}$
with spin $1$ boundary condition on those vertices sitting furthest away
from $x$ in $\Gamma_{(x,n)}$, i.e., those at distance $n$ from 
$\rho$ in $\Gamma$. By the construction of Gibbs measures in Section
\ref{subsect:Potts_infinite}, we have 
\begin{equation}  \label{b's_as_limits}
\lim_{n\rightarrow \infty} b_{x,n} = b_x \, . 
\end{equation}
Imagine now the modified graph $\Gamma^*_{(x,n)}$ obtained from 
$\Gamma_{(x,n)}$ by removing all edges incident to $x$. In other words,
$\Gamma^*_{(x,n)}$ is a disconnected graph with an isolated vertex $x$
together with $k$ connected components isomorphic to
$\Gamma_{(y_1,n)}, \ldots, \Gamma_{(y_k,n)}$. When picking
$X \in \{1, \ldots, q\}^{V_{\Gamma^*_{(x,n)}}}$ according to
$\pi^{\Gamma^*_{(x,n)},1}_{q, \beta}$, the spin configurations on different
connected components obviously become independent. In particular, 
if we only consider the spins $(X(x), X(y_1), \ldots, X(y_k))$, then
we can note that these spins become
independent, with $X(x)$
having distribution $(\frac{1}{q}, \ldots, \frac{1}{q})$, and
$X(y_i)$ having distribution $(\frac{b_{y_i, n}}{b_{y_i, n}+q-1}, 
\frac{1}{b_{y_i, n}+q-1}, \ldots, \frac{1}{b_{y_i, n}+q-1})$. 

If we now reinsert the edges between $x$ and $y_1 , \ldots, y_k$,
thus recovering $\Gamma_{(x,n)}$, then the 
$\pi^{\Gamma_{(x,n)},1}_{q, \beta}$-distribution of
$(X(x), X(y_1), \ldots, X(y_k))$
becomes the same as the corresponding 
$\pi^{\Gamma^*_{(x,n)},1}_{q, \beta}$-distribution
above except that each configuration  
$\xi\in\{1, \ldots, q\}^{\{x, y_1, \ldots, y_k\}}$ is reweighted by a factor
$\exp(2\beta\sum_{i=1}^k I_{\{\xi(y_i)=\xi(x)\}})$. Hence
\[
\pi^{\Gamma_{(x,n)},1}_{q, \beta}((X(x), X(y_1), \ldots, X(y_k))=\xi) 
=\frac{1}{Z}
\prod_{i=1}^k (e^{2\beta I_{\{\xi(y_i)=\xi(x)\}}}
\,b_{y_i,n}^{I_{\{\xi(y_i)=1\}}})
\]
for some normalizing constant $Z$.
By integrating out $X(y_1), \ldots, X(y_k)$, we get 
\[
b_{x,n} = \frac{\prod_{i=1}^k (e^{2\beta} b_{y_i,n} + q - 1)}{
\prod_{i=1}^k (e^{2\beta} + b_{y_i,n} + q - 2)} \, .
\]
Sending $n\rightarrow \infty$ in this expression,
and using (\ref{b's_as_limits}) $k+1$ times
(substituting $x$ with itself and with $y_1, \ldots, y_k$), we obtain
(\ref{eq:wired_recursion}), as desired. $\Cox$

\medskip\noindent
Note that the above proof yields that given $X(x)=1$, the
spins $X(y_1), \ldots, X(y_k)$ become conditionally independent,
with $X(y_i)$ having 
distribution 
\[
\left(\frac{b_{y_i}e^{2\beta}}{b_{y_i}e^{2\beta} + q-1}, 
\frac{1}{b_{y_i}e^{2\beta} + q-1}, \ldots, 
\frac{1}{b_{y_i}e^{2\beta} + q-1}\right) \, .
\]
Likewise, for $l \neq 1$,
conditioning on $X(x)=l$ makes $X(y_1), \ldots, X(y_k)$
conditionally independent with $X(y_i)$ taking value $1$
with probability $\frac{b_{y_i}}{b_{y_i}+ e^{2\beta}+q-2}$,
value $l$ with probability $\frac{e^{2\beta}}{b_{y_i}+ e^{2\beta}+q-2}$,
and other values with probabilities 
$\frac{1}{b_{y_i}+ e^{2\beta}+q-2}$.

By iterating the above argument, we arrive at the following tree-indexed
Markov chain description of the Gibbs measure $\pi^{\Gamma,1}_{q, \beta}$.
\begin{lem}  \label{lem:wired_as_a_Markov_chain}
Suppose that the random spin configuration 
$X\in \{1, \ldots, q\}^{V_\Gamma}$ is obtained as an
inhomogeneous tree-indexed
Markov chain with initial distribution 
\[
\nu=\left(\frac{b_\rho}{b_\rho + q-1}, \frac{1}{b_\rho + q-1}
\ldots, \frac{1}{b_\rho + q-1} \right)
\]
and transition matrices
$P^x = (P^x_{ij})_{i,j \in \{1, \ldots, q\}}$ given by
\[
P^x_{ij} = \left\{ 
\begin{array}{ll}
\frac{b_x e^{2\beta}}{b_x e^{2\beta} + q-1} & \mbox{if } i=j=1 \\
\frac{1}{b_x e^{2\beta} + q-1} & \mbox{if } i=1, \, j \neq 1 \\
\frac{b_x}{b_x+ e^{2\beta}+q-2} & \mbox{if } i\neq 1, j =1 \\
\frac{e^{2\beta}}{b_x+ e^{2\beta}+q-2} & \mbox{if } i=j \neq 1 \\
\frac{1}{b_x+ e^{2\beta}+q-2} & \mbox{otherwise.} 
\end{array} \right.
\]
Then $X$ has distribution $\pi^{\Gamma,1}_{q, \beta}$.
\end{lem}
A crucial difference now compared to the Gibbs measure 
$\pi^{\Gamma,0}_{q, \beta}$ with free boundary condition
is that if any $b_x \neq 1$, then there is not enough state-symmetry
in the tree-indexed Markov chain in Lemma \ref{lem:wired_as_a_Markov_chain}
to make the corresponding fuzzy Potts model a tree-indexed Markov chain.
This will soon become clear. 

A key lemma for proving nonquasilocality in the fuzzy Potts model is
the following. 
\begin{lem}  \label{lem:sibling}
If $\pi^{\Gamma,1}_{q, \beta} \neq \pi^{\Gamma,0}_{q, \beta}$, then
there exist two siblings $y_1, y_2 \in V_\Gamma$ such that
$b_{y_i}>1$ for both $i=1$ and $i=2$. 
\end{lem}
{\bf Proof:}
It follows from the assumption 
$\pi^{\Gamma,1}_{q, \beta} \neq \pi^{\Gamma,0}_{q, \beta}$ using 
(\ref{eq:positive_magnetization}) that $a_\rho>\frac{1}{q}$, so that
\begin{equation}  \label{eq:magnetization_at_rho}
b_\rho>1 \, .
\end{equation}
Furthermore, (\ref{eq:positive_magnetization}) and
(\ref{eq:zero_magnetization}) imply that $a_x \geq \frac{1}{q}$ for
all $x\in V_\Gamma$, whence $b_x \geq 1$ for
all $x\in V_\Gamma$. Note also that $1$ is a fixed point
of the recursion (\ref{eq:wired_recursion}), in the sense that
if all children $y_1, \ldots, y_k$ satisfy $b_{y_i}=1$, then $b_x=1$. 

Hence, $\rho$ must have at least one child $x$
with $b_x>1$. By iterating this argument we see that for any $n$,
it must have at least one descendant $x$ at distance $n$ such that
$b_x>1$. Fix $n$ and such a vertex $x$ with $b_x>1$
at distance $n$ from $\rho$. Write $(z_0, z_1, \ldots, z_n)$ for the
vertices on the self-avoiding path from $x$ to $\rho$ (so that in 
particular $z_0=x$ and $z_n = \rho$). Next, note that
the recursion (\ref{eq:wired_recursion}) has the property that if
one of the children $y_i$ has $b_{y_1}>1$, then $b_x>0$ as well. Since
$b_{z_0}>1$ it follows that $b_{z_i}>1$ for $i=1, \ldots, n$. 

Suppose now for contradiction that the assertion of the lemma is false,
i.e., that there are no two siblings $y_1, y_2 \in V_\Gamma$ for which
$b_{y_1}>1$ and $b_{y_2}>1$. Then none of the vertices $z_0, \ldots,
z_{n-1}$ has a sibling $y$ with $b_y>1$. The recursion 
(\ref{eq:wired_recursion}) along the path $(z_0, z_1, \ldots, z_n)$
then turns into a simple one-dimensional dynamical system on the space
$[1, \infty)$ given by $b_{z_{i+1}}= f(b_{z_i})$ where
\[
f(b) = \frac{e^{2\beta}b+q-1}{e^{2\beta}+b+q-2}  \, .
\]
This dynamical system is contractive with a unique fixed point at
$b=1$, so that -- if we just keep iterating beyond the $n$'th iteration --
for any initial value $b_{z_0}\in [1, \infty)$ we obtain
\begin{equation}  \label{eq:limit_of_dynamical_system}
\lim_{n \rightarrow \infty} b_{z_n}=1 \, .
\end{equation}
Since $f$ is increasing and bounded by $e^{2\beta}$, we get that
$b_{z_1}$ is bounded by $e^{2\beta}$ and, therefore, that the
convergence in (\ref{eq:limit_of_dynamical_system}) is in fact uniform
in the initial value $b_{z_0}$.
Thus we can, for any $\eps>0$, find an $n$ which guarantees that
$b_{z_n}<1+ \eps$. Thus, $b_\rho <1+ \eps$ for any $\eps>0$, whence
$b_\rho=1$. But this contradicts (\ref{eq:magnetization_at_rho}), so
the proof is complete. $\Cox$

\medskip\noindent
{\bf Proof of Proposition \ref{prop:wired}:}
By Lemma \ref{lem:sibling}, $\Gamma$ has at least one vertex which
has (at least) two children $y_1$ and $y_2$ that both have 
$b_{y_i}>1$. The choice of root $\rho$ for the tree does not influence
the Gibbs measure $\pi^{\Gamma,1}_{q, \beta}$, and therefore we may assume
that $\rho$ has two such children $y_1$ and $y_2$. 
We shall for simplicity first prove the proposition under the assumption
that 
\begin{equation}  \label{eq:preliminary_assumption}
\mbox{$\rho$ has no other children,}
\end{equation}
and in the end show how to remove this assumption. 

We shall have a look at the conditional distribution of the fuzzy spin
$Y(\rho)$ at the root, given that its neighbors (i.e., its children)
take value 
\begin{equation}  \label{eq:all_1's}
Y(y_1)= Y(y_2)=1 \, .
\end{equation}
By summing over all 
$X\in \{1, \ldots, q\}^{\{\rho, y_1, y_2\}}$ such that
(\ref{eq:all_1's}) holds, and using Lemma \ref{lem:wired_as_a_Markov_chain},
we obtain
\begin{eqnarray}  \label{eq:crucial_ratio_of_probabilities}
\lefteqn{\frac{\P(Y(\rho)=1 \, | \, Y(y_1) = Y(y_2)=1)}{
\P(Y(\rho) \neq 1 \, | \, Y(y_1) = Y(y_2)=1)}} \\  \nonumber
& = & \frac{\frac{b_\rho}{b_\rho+q-1}\prod_{i=1}^2 
\frac{b_{y_i}e^{2\beta}+r_1-1}{{b_{y_i}e^{2\beta}+q-1}} +
\frac{r_1-1}{b_\rho+q-1}\prod_{i=1}^2 
\frac{b_{y_i}+e^{2\beta}+r_1-2}{{b_{y_i}+e^{2\beta}+q-2}}}
{\frac{q-r_1}{b_\rho+q-1} \prod_{i=1}^2 
\frac{b_{y_i}+r_1-1}{{b_{y_i}+e^{2\beta}+q-2}}}  \\  \nonumber
& = & 
\frac{\prod_{i=1}^2 (b_{y_i}e^{2\beta} + r_1-1)+ (r_1-1)\prod_{i=1}^2
 (b_{y_i}+e^{2\beta} + r_1-2)}{(q-r_1) \prod_{i=1}^2 (b_{y_i} +r_1-1)}
\end{eqnarray}
where in the last line we have used (\ref{eq:wired_recursion}) to express
$b_\rho$ in terms of the $b_{y_i}$'s. 

Now pick an $n$, and consider conditioning further on some 
$\eta_n \in \{1, \ldots, s\}^{V_{\Gamma_{n+1}} \setminus \{\rho\}}$ 
such that $\eta_n(y_1) = \eta_n(y_2)=1$. 
The conditional probability $\P(Y(\rho)=1 \, | \, Y(y_1) = Y(y_2)=1)$
is a convex combination of terms
$\P(Y(\rho)=1 \, | \, Y(V_{\Gamma_{n+1}} \setminus \{\rho\}) = \eta_n)$
for such $\eta_n$'s. We can therefore find a particular 
$\eta_n\in\{1, \ldots, s\}^{\Lambda_{n+1} \setminus \{\rho\}}$ such that
\begin{eqnarray}
\nonumber
\lefteqn{
\frac{\P(Y(\rho)=1 \, | \, Y(V_{\Gamma_{n+1}} \setminus \{\rho\}) = \eta_n)}{
 \P(Y(\rho) \neq 1 \, | \, Y(V_{\Gamma_{n+1}} \setminus \{\rho\}) 
= \eta_n)}} \\
& \geq &
\frac{\prod_{i=1}^2 (b_{y_i}e^{2\beta} + r_1-1)+ (r_1-1)\prod_{i=1}^2
 (b_{y_i}+e^{2\beta} + r_1-2)}{(q-r_1) \prod_{i=1}^2 (b_{y_i} +r_1-1)}
\label{eq:wired_conditioning_on_eta} \, .
\end{eqnarray}
Fix such an $\eta_n$. Next, construct another configuration 
$\eta'_n\in \{1, \ldots, s\}^{V_{\Gamma_{n+1}} \setminus \{\rho\}}$ by taking
\[
\eta'_n(x) = \left\{
\begin{array}{ll}
\eta_n(x) & \mbox{for } x \in  V_{\Gamma_n} \setminus \{\rho\}  \\
(\eta_n(\parent(x))+1) \, {\rm mod} \, s & \mbox{for } x \in V_{\Gamma_{n+1}}
\setminus V_{\Gamma_n} \, .
\end{array} \right. 
\]
The crucial aspects of this choice of $\eta'_n$ is that (a) $\eta_n=\eta'_n$ on
$V_{\Gamma_n}$ and (b) each $x$ in the remotest layer 
$V_{\Gamma_{n+1}}\setminus V_{\Gamma_n}$ of $\Gamma_{n+1}$ has a 
fuzzy spin value
which is different from its parent. It is readily checked that property (b)
implies that the conditional distribution of
$Y(V_{\Gamma_{n-1}})$ given $Y(V_{\Gamma_{n+1}}\setminus V_{\Gamma_{n-1}})=
\eta'_n(V_{\Gamma_{n+1}}\setminus V_{\Gamma_{n-1}})$ becomes the same
as if the underlying Gibbs measure had been not $\pi^{\Gamma,1}_{q, \beta}$
but rather the finite-volume Gibbs measure $\pi^{\Gamma_{n+1}}_{q, \beta}$
(cf.\ \cite[Lem.\ 9.2]{H3}). Hence the conditional distribution
of $Y(\rho)$ given that $Y(V_{\Gamma_{n+1}} \setminus \{\rho\}) = \eta'_n)$ 
can be calculated from the tree-indexed Markov chain corresponding to
free boundary condition, i.e., the one defined in 
(\ref{Markov_chain_for_Y_with_free_bc}). We get 
\begin{equation}  \label{eq:free_conditioning_on_eta} 
\frac{\P(Y(\rho)=1 \, | \, Y(V_{\Gamma_{n+1}} \setminus \{\rho\}) = \eta'_n)}{
 \P(Y(\rho) \neq 1 \, | \, Y(V_{\Gamma_{n+1}} \setminus \{\rho\}) = \eta'_n)} 
= 
\frac{(e^{2\beta}+r_1-1)^2}{(q-r_1)r_1}
\end{equation}
Note that the right-hand sides of 
(\ref{eq:wired_conditioning_on_eta}) and (\ref{eq:free_conditioning_on_eta})
do not depend on $n$. 
We now make the following crucial claim.
\begin{quote}
{\bf Claim:} If $b_{y_1}>1$ and $b_{y_2}>1$, then
the right-hand side of (\ref{eq:wired_conditioning_on_eta}) is strictly
greater than the right-hand side of (\ref{eq:free_conditioning_on_eta}).
\end{quote}
To prove the claim, define
\[
a = \frac{b_{y_1}b_{y_2}+r_1 - 1}{(b_{y_1}+r_1 - 1)(b_{y_2} + r_1 - 1)}
\]
and note that $a$ can be rewritten as
\begin{eqnarray*}
a & = & \frac{b_{y_1}b_{y_2}+r_1 - 1}{(b_{y_1}+r_1 - 1)(b_{y_2} + r_1 - 1)} \\
& = & \frac{1}{r_1}\frac{r_1}{(b_{y_1}+r_1-1)} +
\frac{b_{y_2}}{(b_{y_2}+r_1 - 1)} \frac{(b_{y_1}-1)}{(b_{y_1}+r_1-1)} \, .
\end{eqnarray*}
Assuming that $b_{y_1}>1$ and $b_{y_2}>1$, we get that 
$\frac{b_{y_2}}{b_{y_2}+r_1 - 1} > \frac{1}{r_1}$ and that
$\frac{b_{y_1}-1}{b_{y_1}+r_1-1}>0$, whence
\begin{eqnarray} \nonumber
a & > & \frac{1}{r_1}\frac{r_1}{(b_{y_1}+r_1-1)} +
\frac{1}{r_1} \frac{(b_{y_1}-1)}{(b_{y_1}+r_1-1)} \\
& = & \frac{1}{r_1} \, .
\label{eq:inequality_for_a}
\end{eqnarray}
Next, an elementary but tedious calculation shows that the right-hand side
of (\ref{eq:wired_conditioning_on_eta}) can be rewritten as
\begin{equation}  \label{eq:rewritten_wired}
\frac{a(e^{4\beta}+r_1-1) + (1-a) (2e^{2\beta} + r_1 - 2)}{q - r_1} \, .
\end{equation} 
Analogously, the right-hand side of (\ref{eq:free_conditioning_on_eta})
can be rewritten as 
\begin{equation}  \label{eq:rewritten_free}
\frac{\frac{1}{r_1}(e^{4\beta}+r_1-1) + (1-\frac{1}{r_1}) 
(2e^{2\beta} + r_1 - 2)}{q - r_1} \, .
\end{equation} 
Now, using (\ref{eq:inequality_for_a}) and the observation that
\[
e^{4\beta}+r_1-1 > 2e^{2\beta} + r_1 - 2 \, ,
\]
we get that the expression in (\ref{eq:rewritten_wired}) is strictly greater
than that in (\ref{eq:rewritten_free}), and the claim is proved. 

Hence the difference between the left-hand sides of
(\ref{eq:wired_conditioning_on_eta}) and (\ref{eq:free_conditioning_on_eta})
is bounded away from $0$ uniformly in $n$.
The denominators of the left-hand sides are bounded away from $0$ uniformly
in $n$ due to uniform nonnullness of the fuzzy Potts model
(see Section \ref{sect:quasilocal}). Hence
\[
\P(Y(\rho)=1 \, | \, Y(V_{\Gamma_{n+1}} \setminus \{\rho\}) = \eta_n)
- \P(Y(\rho)=1 \, | \, Y(V_{\Gamma_{n+1}} \setminus \{\rho\}) = \eta'_n)
\]
is bounded away from $0$ uniformly on $n$. By plugging in 
these $\eta_n$ and $\eta'_n$ in (\ref{eq:def_quasilocality}), we get,
since $\eta_n=\eta'_n$ on $V_{\Gamma_n}$, that quasilocality of $Y$
fails. This proves the proposition modulo the assumption 
(\ref{eq:preliminary_assumption}). 

It remains to remove the assumption (\ref{eq:preliminary_assumption}). 
To do this, suppose that $\rho$ has $k-2$ additional children 
$y_3, \ldots, y_k$. We can then extend the configurations 
$\eta$ and $\eta'$ that we
condition on above, to $y_3, \ldots, y_k$ and their descendants, as follows.
We insist that $\eta$ and $\eta'$ 
that they take value $1$ at $y_3, \ldots, y_k$, and that they take some value 
other than $1$ at all children of $y_3, \ldots, y_k$ (they may
otherwise be arbitrary on the further descendants of
$y_3, \ldots, y_k$). 
Easy modifications of the calculations above show that 
(\ref{eq:wired_conditioning_on_eta}) and (\ref{eq:free_conditioning_on_eta})
hold as before, with the modification that both right-hand sides are
multiplied by 
\[
\left( \frac{e^{2\beta} + r_1 - 1}{r_1} \right)^{k-2} \, .
\]
Since this factor is the same in 
(\ref{eq:wired_conditioning_on_eta}) and (\ref{eq:free_conditioning_on_eta}),
the rest of the proof goes through as before. 
$\Cox$

\medskip\noindent
{\bf Remark:} Since the event conditioned on in 
(\ref{eq:crucial_ratio_of_probabilities}) has positive measure, it is
easy to extract from the above proof that the set of discontinuities
of 
the conditional probability
$\P(Y(\rho)=1 | Y(V_\Gamma\setminus \{\rho\}) = \eta)$ as a function of
$\rho$, has positive measure under $\mu^{G,1}_{q, \beta, (r_1, \ldots, r_s)}$.
Hence, so-called almost sure quasilocality and almost sure Gibbsianness
fails in general 
for the fuzzy Potts model on trees, in contrast to the $\Z^d$ case 
(see Maes and Vande Velde \cite{MVV}) and the mean-field case
(Theorem \ref{thm:main_result_on_complete_graph} (iv)). This 
contrast between the fuzzy Potts model on $\Z^d$ and on trees is analogous
to the corresponding almost sure Gibbsianness issue for the random-cluster
model; see \cite{H96}. 

\subsection{Discussion}

What concrete information can we extract from 
Theorem \ref{thm:main_result_on_trees}? Let $\beta_c=\beta_c(\Gamma,q)$
denote, as in Section \ref{subsect:Potts_infinite}, the critical value
for the $q$-state Potts model on the tree $\Gamma$. For $q\geq 3$, 
we then have
from Theorem \ref{thm:main_result_on_trees} that $\beta<\beta_c$ implies
that any corresponding fuzzy Potts measure is Gibbsian, while
$\beta>\beta_c$ yields existence of corresponding fuzzy Potts measures
that are non-Gibbsian. 

It remains to specify the critical value $\beta_c(\Gamma,q)$. If we know the 
critical value $p_c(\Gamma,q)$ of the corresponding random-cluster
model, then we can calculate $\beta_c= -\frac{1}{2} \log (1-p_c)$
(see, e.g., \cite{GHM}). For the case when $\Gamma$ is a regular tree,
the critical value $p_c(\Gamma,q)$ can be characterized
in terms of the solutions of a certain algebraic equation given in
\cite[p.\ 235]{H96a}.

For general trees the situation is more complicated. For a variety of
stochastic models on trees, critical values can be calculated
in terms of a natural quantity known as the branching number of
the tree, denoted ${\rm br}(\Gamma)$; see for instance \cite{P}. 
Lyons \cite{L} calculated $\beta_c(\Gamma,q)$ in terms of ${\rm br}(\Gamma)$
for the case $q=2$. In contrast, and perhaps somewhat surprisingly, 
the critical values $\beta_c(\Gamma,q)$ for
larger $q$ do {\em not} admit a characterization in terms of
${\rm br}(\Gamma)$; this was shown by Pemantle
and Steif \cite{PS}. Bounds for $\beta_c(\Gamma,q)$ that only depend on 
${\rm br}(\Gamma)$ and on $q$ can, however, be obtained using the
standard comparison techniques for the random-cluster model reviewed
in \cite{GHM}.

\section{The fuzzy Potts model on complete graphs}  \label{sect:mean-field}

In this section we treat the case of complete graphs. 
We start with precise definitions of the model and 
a detailed explanation of the limiting process for the conditional 
probabilities that was sketched in the introduction. 
The proofs are 
essentially self-contained but use some standard knowledge 
(whose main reference is Ellis and Wang \cite{ElWa90}) on the 
infinite volume limit of the empirical distribution of the order 
parameter in the mean-field Potts model.

\subsection{Mean-field Potts in finite volume $N$}

For a positive integer $q$,   the
Gibbs measure $\pi^N_{q,\beta}$ for the $q$-state Potts model on the complete 
graph with $N$ vertices 
at inverse temperature $\beta\geq 0$, is
the probability measure  on $\{1, \ldots, q\}^N$ 
which to each $\xi\in \{1, \ldots, q\}^N$ assigns probability
\begin{equation} \label{eq:Potts-mf}
\pi^N_{q,\beta}(\xi)=\frac{1}{Z^N_{q, \beta}}\exp\left(
\frac{\beta}{N} \sum_{1\leq x\neq y \leq N} I_{\{\xi(x) = \xi(y)\}} \right) \, .
\end{equation}
Here $Z^N_{q, \beta}$ is the
normalizing constant. 
Note that this definition  slightly deviates from the definition 
(\ref{eq:Potts}) by the factor $1/N$ appearing in the exponential. 
Such a convention is appropriate because, clearly,  the interaction 
must be chosen depending on the size of the graph in a mean-field model. 
This definition of the finite volume Gibbs-measures 
is standard in the literature; see e.g.\ \cite{ElWa90}.


\subsection{Mean-field fuzzy Potts in finite volume $N$}

The mean-field fuzzy Potts measure in finite volume $N$ is then defined in the same way as it is defined 
on every graph. 
To be explicit,  fix $q$, $\beta$ and the spin-partition 
$(r_1, \ldots, r_s)$ as above. 
 Let $X$ be the
$\{1, \ldots, q\}^N$-valued random object distributed according to the mean field finite volume Gibbs
measure $\pi^N_{q, \beta}$. Take $Y$ to be the 
$\{1\ldots, s\}^N$-valued random object obtained from $X$ by the site-wise application 
of the spin-partitioning as in
(\ref{eq:def_of_fuzzy_Potts}). 
Then $\mu^N_{q, \beta, (r_1, \ldots, r_s)}$ is the probability
measure on $\{1, \ldots, s\}^N$ which describes the distribution of $Y$.

\subsection{Gibbsianness vs.\ non-Gibbsianness for mean-field models: \\
continuity vs.\ discontinuity of limiting conditional probabilities}

We start with some general remarks about mean-field models 
to explain the appropriate analogue of non-Gibbsianness in more detail than 
we did in the introduction. 
To begin with, the following lemma makes 
explicit that we can always describe the single-site conditional 
probabilities of the finite volume Gibbs measures of a mean-field model  
in terms of a {\it single-site kernel} from the empirical distribution 
vector of the conditioning to the single-site state space. 
It is the infinite volume limit of this kernel that shall then be considered 
in the analysis of the model. 

So, suppose that $S$ is a finite set (local spin space) 
and for any $N$ we are given an exchangeable (that is permutation-invariant) measure 
$\mu^N$ on $S^N$. This permutation invariance
is certainly true for the mean-field Potts model. Moreover it carries over trivially  
to the fuzzy Potts model. This is clear since the distribution of the latter is simply
obtained by an application of the same map to the spin variable at each site. 

In a general context, 
denote by $\PP=\{ (p_i)_{i\in S}, 0\leq p_i\leq 1, \sum_{i\in S}p_i =1  \}$ the space of probability 
vectors on the set $S$.
We use the obvious short notation $x^c=\{1,\dots, N\}\backslash \{x\}$.

\begin{lem}  \label{lem:kernel}
For each $N$ there is a probability kernel $Q^N:S\times\PP\rightarrow [0,1]$
from $\PP$ to the single-site state space $S$ such that the single-site 
conditional expectations at any site $x$ can be written in the form 
\begin{equation} \label{eq:0.1}
\mu^N\bigl(X(x)=i\bigl | X(x^c)=\eta \bigr)=Q^N
\bigl(i\bigl | (n_j)_{j\in S} \bigr)
\, .
\end{equation}
Here $n_j=\frac{1}{N-1}\#\bigl(1\leq y\leq N, y\neq x, \eta(x)=j\bigr)$ 
is the fraction of sites for which the spin-values of the conditioning 
are in the state $j\in S$.  
\end{lem}
 
\noindent {\bf Proof:}
By exchangeability it is clear that the right hand side 
of  (\ref{eq:0.1}) depends on
the sets $\bigl\{1\leq y\leq N, y\neq x, \eta(y)=j\bigr\}$, for all $j\in S$,  
only through their size. Equivalently we may express this dependence in terms 
of the empirical distribution $(n_j)_{j\in S}$.
$\Cox$

\medskip\noindent
In turn, the knowledge of the kernel $Q^N$ uniquely determines 
the measure $\mu^N$. This is clear since the knowledge of all one-site conditional 
probabilities of finitely many random variables uniquely determines the joint distribution. 
So we may as well consider the $Q^N$«s as the basic objects and regard them 
as the starting point of the definition of a mean field model. 
This is of course only meaningful if the $Q^N$«s are related to each other in 
a meaningful way. 

Let us turn now to the concrete case of the mean-field Potts model to point 
out two very simple observations that shall serve as a motivation 
of our further investigation. In this case we have directly from 
the definition 
(\ref{eq:Potts-mf}) the explicit formula 
\begin{equation} \label{eq:0.2}
Q^N_{q,\b}
\bigl(i\bigl | (n_j)_{1\leq j\leq q} \bigr)=\frac{\exp\bigl( \b (1-\frac{1}{N})n_i  \bigr)}{
\sum_{j=1}^q \exp\bigl( \b (1-\frac{1}{N})n_j  \bigr)
}
\, .
\end{equation}
We note the following. 
\begin{description}
\item{\bf (i)} $Q^N_{q,\b}$ converges to 
$Q^\infty_{q,\b}  =\frac{\exp\bigl( \b  n_i  \bigr)}{
\sum_{j=1}^q \exp\bigl( \b n_j  \bigr)}$ when $N$ tends to infinity. 
Indeed, the trivial $1/N$-factor appearing in (\ref{eq:0.2}) could of course 
even be removed by a harmless 
redefinition of the model that would lead to the same 
infinite volume behavior of the Gibbs measures, making all $Q^N_{q,\b}$ identical. 
\item {\bf (ii)} The limiting kernel 
$Q^\infty_{q,\b}$ is a {\it continuous function} of the probability 
vector $(n_j)_{1\leq j\leq q}$, as a function on $\R^q$. 
\end{description}
The existence of the infinite volume limit (i) is a minimal ingredient 
for the definition of a mean-field model. 
Assuming this 
we can talk about limiting or ``infinite volume'' conditional probabities. 
Then, {\it continuous dependence of the limiting conditional probability} as it is stated 
in (ii) is the  
obvious analogue to the continuous dependence of the conditional 
expectation of a lattice model on the conditioning with respect to 
product topology. 

So, properties (i) and (ii) are the analogues of a proper 
Gibbsian structure for mean-field models. 
``Non-Gibbsianness'' may then manifest 
itself by the failure of (ii) at certain 
points of discontinuity.  The reader may find a number of examples 
of this in \cite{K3}. 
After these introductory remarks 
we will show in the following that discontinuities in fact occur
for the mean-field fuzzy Potts model,  for certain values of the parameters,  
and discuss them in detail.

\subsection{Conditional probabilities for fuzzy Potts in finite volume}

Let us use the following notation for the single-site 
probability kernel that describes the conditional probabilities of 
the fuzzy model.

\begin{equation} \label{eq:0.3}
\mu^N_{q, \beta, (r_1, \ldots, r_s)}\Bigl(Y(x)=k\bigl | Y(x^c)=\eta \Bigr)=:Q^N_{q, \beta, (r_1, \ldots, r_s)}
\bigl(k\bigl | (n_l)_{1\leq l\leq s} \bigr)
\, .
\end{equation}
where $n_l=\frac{1}{N-1}\#\bigl(1\leq y\leq N, y\neq x, \eta(x)=l\bigr)$, for $l=1,\dots,s$ 
is the empirical distribution of  fuzzy spin-values in the conditioning. 

Now, it is not difficult to derive an explicit expression in 
terms of expectations with respect to 
ordinary mean-field Potts measures, having the number 
of states given by the sizes of the classes $r_l$.  
Clearly, the infinite volume analysis relies on this result.

\begin{prop}  \label{prop:1}For each finite $N$ we have the representation 
\begin{equation} \label{eq:prop_1}
Q^N_{q, \beta, (r_1, \ldots, r_s)}
\bigl(k\bigl | (n_l)_{1\leq l\leq s} \bigr)
=\frac{r_{k} \,A(\b_{k},r_{k},N_{k})
}{\sum_{l=1}^{s}r_{l} \,A(\b_{l},r_{l},N_{l})
}
\end{equation}
where
\[
A(\tilde\b,r,M) \equiv  \pi^M_{r,\tilde\beta}\Bigl( 
\exp\Bigl( \frac{\tilde\b}{M} \sum_{x=1}^{M} I_{X(x)=1}
\Bigr)
\Bigr) \, ,
\]
\[
N_{k} = (N-1)n_k  \, ,
\]
and 
\[
\b_{k} = \frac{\b N_{k}}{N} =\b \Bigl(1-\frac{1}{N}\Bigr)n_k \, .
\]
\end{prop}

\noindent{\bf Remark:} In particular we have $A(\tilde\b,r=1,N)=e^{\tilde\b}$. 
From this we see immediately that the case of the original Potts model 
is recovered by setting all $r_l$ equal to one. 

\medskip\noindent 
{\bf Proof of Proposition \ref{prop:1}:}
To compute the left hand side of (\ref{eq:prop_1})
we may choose $x=1$ and write 
\begin{eqnarray*}
\lefteqn{\textstyle \mu^N_{q, \beta, (r_1, \ldots, r_s)}\Bigl(
Y(1)=k\Bigl |Y([2,N])=\eta([2,N])\Bigl)} \\
&=& {\textstyle \frac{1}{\Norm(\eta([2,N])) }
\sum_{\xi(1)\mapsto k} 
\sum_{\xi([2,N])\mapsto \eta([2,N])} 
\pi^N_{q,\beta}\Bigl(\xi(1),\xi([2,N])\Bigr) \, .}
\end{eqnarray*}
Here we are summing over Potts configurations $\xi$ that are mapped 
to the fuzzy Potts configuration $(k,\eta)$ by means of the definition 
of the fuzzy model given in (\ref{eq:def_of_fuzzy_Potts}). 
The normalization has to be 
chosen such that 
summing over $k=1,\dots,s$ yields one, for each fixed $\eta([2,N])$. 
The partition function appearing in the Gibbs-average on the right
hand side only gives 
a constant that can be absorbed in the normalization, and so we need only 
consider  
\begin{eqnarray*}
\lefteqn{\textstyle \sum_{\xi(1)\mapsto k} 
\sum_{\xi([2,N])\mapsto \eta([2,N])}
\exp\Bigl( 
\frac{\beta}{N} \sum_{1\leq x\neq y \leq N} I_{\{\xi(x) = \xi(y)\}}
\Bigr)}\hspace{15mm} \\
&= & {\textstyle \sum_{\xi(1)\mapsto k} 
\sum_{\xi([2,N])\mapsto \eta([2,N])}\exp\Bigl( 
 \frac{\b}{N}\sum_{2\leq y\leq N} I_{\{\xi(1)=\xi(y)\}}
\Bigr)} \\
& & {\textstyle 
\times
\exp\Bigl( 
 \frac{\b}{N}\sum_{2\leq x\neq y \leq N} I_{\{\xi(x) = \xi(y)\}}
\Bigr) \, . }
\end{eqnarray*}

For fixed $\eta([2,N])$ we denote $\L_{l}:= \#\bigl\{ 
x\in \{2,\dots,N\}:\eta(x)=l \bigr\}$. 
Then the sum in the last exponential decomposes over these sets, 
and we can rewrite the right hand side of the last equation in the form  
\begin{eqnarray*}
&
\sum_{\xi(1)\mapsto k} 
\sum_{\xi([2,N])\mapsto \eta([2,N])}
\exp\Bigl(\frac{\b_{k}}{N_{k}} \sum_{z\in \L_{k}} I_{\{\xi(z)=\xi(1)\}}\Bigr)\cr
&\quad\times \prod_{l=1}^s \exp\Bigl(\frac{\b_{l}} {N_{l}}
\sum_{x<y, x,y\in \L_l} 
I_{\{\xi(x) = \xi(y)\}}\Bigr) \, .
\end{eqnarray*}
Next we divide 
the last line by the product of partition functions which is obtained 
by omitting the first exponential and the first sum. 
This only yields another $\eta([2,N])$-dependent constant.
Using cancellations for the terms 
with $l\neq k$ we see in this way that 
\begin{eqnarray*}
\lefteqn{\mu^N_{q, \beta, 
(r_1, \ldots, r_s)}\Bigl(Y(1)=k\Bigl |Y([2,N])=\eta([2,N])\Bigl)} 
\hspace{10 mm} \\
&=& \frac{1}{\Norm(\eta([2,N])) }
\sum_{\xi(1)\mapsto k}  \pi^{N_k}_{k,\beta_k}\Bigl( 
\exp\Bigl( \frac{\b_k}{N_k} \sum_{z=1}^{N_k} I_{\{X(z)=1\}}
\Bigr)
\Bigr) \, ,
\end{eqnarray*}
which concludes the proof.
$\Cox$

\subsection{Continuity vs.\ 
discontinuity of limiting conditional probabilities for fuzzy Potts}

In this subsection we will derive an explicit formula for the limiting conditional 
probabilities of the fuzzy model. From this parts (i), (ii), (iii) of Theorem 
 \ref{thm:main_result_on_complete_graph} follow. 

We can build on well-known results 
about the limiting behavior of the 
empirical distribution of the mean-field Potts model. 
The main point is that it exhibits a first-order 
phase transition at a finite inverse critical 
temperature $\b_c(q)$, for all $q\geq 3$. For the special 
case $q=2$ (Ising model) there is only 
a second order phase transition. 
The following pieces of information 
about the mean-field Potts model can be 
found in \cite[Thms 2.1. and 2.3]{ElWa90}.
The reader should focus at first on the case 
$\b\neq \b_c(q)$, i.e. off the critical temperature. 

\begin{thm}  {\bf (Ellis, Wang)}\label{thm:Empirical measures}
Assume that $q\geq 3$, and
suppose that 
$\b\neq \b_c(q):= \frac{2(q-1)}{q-2}\log(q-1)$.  Then we have the weak limit
\begin{eqnarray} \label{eq:emp-meas}
\lefteqn{\lim_{N\uparrow\infty}
\pi^N_{q,\beta}\Bigl(
\frac{1}{N}\sum_{x=1}^N (I_{\{X(x)=1\}},\dots, I_{\{X(x)=q\}})
\in \cdot \Bigr)} \\ 
&= & \left\{ 
\begin{array}{ll}
\d_{\frac{1}{q}(1,1,\dots,1)}, & \mbox{if } \b<\b_c(q) \\
\frac{1}{q} \sum_{\nu=1}^q\d_{u(\b,q) \, e_{\nu}+\frac{1-u(\b,q)}{q}(1,1,\dots,1)}, & \mbox{if } \b>\b_c(q)\\
\l_0(q)\d_{\frac{1}{q}(1,1,\dots,1)}
+\frac{1-\l_0(q)}{q} \sum_{\nu=1}^q\d_{u(\b_c(q),q) \, e_{\nu}+\frac{1-u(\b_c(q),q)}{q}(1,1,\dots,1)}
& \mbox{if } \b=\b_c(q) \, ,
\end{array} \right.
\nonumber
\end{eqnarray}
where $e_{i}$ is the unit vector in the $i$'th coordinate 
direction of $\R^q$.

The quantity $u(\b,q)$ is well defined for $\b\geq \b_c(q)$. 
It is the largest solution 
of the mean field equation 
\begin{eqnarray} \label{eq:meanfield-eq}u=\frac{1-e^{-\b u}}{
1+(q-1)e^{-\b u}
}
\end{eqnarray}
and obeys the following properties:  
It is strictly increasing in $\b$,  and we have  $u(q,\b_c(q))=\frac{q-2}{q-1}$.
The constant appearing at the critical point obeys the strict inequality $0<\l_0(q)<1$. 
\end{thm}

Some comments are in order:  Obviously, $u(\b,q)$ plays the role of an 
order parameter. 
Now, for 
$\b>\b_c(q)$ the system is in a symmetric linear combination of $\nu$-like states.  
The limiting empirical distribution becomes the equidistribution on the possible 
spin values for $\b<\b_c(q)$. It jumps at the critical point for $q\geq 3$.
At the critical point itself there is a non-trivial 
linear combination between both types of measures. 

To feel comfortable with the mean-field equation (\ref{eq:meanfield-eq}) the reader may 
note that it is obtained from the  equations $n_i=\frac{\exp\bigl( \b n_i  \bigr)}{
\sum_{j=1}^q \exp\bigl( \b n_j  \bigr)}$ for $i=1,\dots,q$ with the following ansatz:
Denote by $i$ the index with the largest $n_j$. Assume that $n_j$ is independent 
of $j$, for $j\neq i$, and put $u=n_i-n_j$ for some $j\neq i$. 

Let us mention that the results of Theorem \ref{thm:Empirical measures}
can be obtained by a Gaussian transformation and saddle point estimates
on the resulting integrals (all of which is omitted here).  
At the critical point a little care is needed: To obtain the proper value of the constant $\l_0(q)$ 
a Gaussian approximation  around the minima 
and estimates showing positive curvature are needed. 

The well-known case of the mean field 
Ising model $q=2$ can be recovered from the theorem 
by taking the formal limit  $q\downarrow 2$ in the explicit formula for $\b_c(q)$
and noting that $u(q,\b_c(q))=0$. So (\ref{eq:emp-meas}) describes 
a second order transition in that case.




The following explicit formula for the limiting conditional 
probabilities of the fuzzy model now follows easily from 
our finite volume representation 
of the conditional probabilities given in Proposition \ref{prop:1}
and the known limiting statement of Theorem \ref{thm:Empirical measures}.

\begin{thm}\label{thm:5.4} We have
\[
\lim_{N\uparrow \infty}Q^N_{q, \beta, (r_1, \ldots, r_s)}
\bigl(k\bigl | (n_l)_{1\leq l\leq s} \bigr)
= \frac{C(\b n_{k},r_{k})
}{\sum_{l=1}^{s} C(\b n_{l},r_{l})
}
\]
whenever $n_{k}\neq \b_{c}(r_{k})/\b$ 
for all $k$ with $r_{k}\geq 3$.
Here 
\[
C(\tilde\b,r)=\exp\Bigl( \frac{\tilde\b}{r} \Bigr)\times \left\{ 
\begin{array}{ll}
 r,& \hbox{if}\,\,  \tilde\b<\b_c(r)
\cr
 \exp\Bigl( \frac{\tilde\b(r-1)u(\tilde\b,r)}{r}\Bigr)
+ (r-1)\exp\Bigl( 
-\frac{\tilde\b \,u(\tilde\b,r)}{r}\Bigr)
,& \hbox{if} \,\, \tilde\b>\b_c(r) \, . 
\end{array}\right.
\]
\end{thm} 


\noindent{\bf Proof of Theorem \ref{thm:5.4}:} 
Let $\tilde \b\neq\b_c(q)$. By  Theorem \ref{thm:Empirical measures}  we have 
$\lim_{M\uparrow\infty}r A(\tilde \b,r,M)=C(\tilde \b,r)$.
$\Cox$ 

\medskip\noindent

\noindent{\bf Remark:} Obviously this gives the right answer for $\b=0$ 
or in the case of the original Potts model (letting all $r_l$ be equal to one). 
We see however that the limiting form of 
the conditional expectations has a nontrivial form in general.  
This expression has jumps for $n_{l}=\b_{c}(r_{l})/\b$ 
whenever $r_{l}\geq 3$. (For matters of simplicity we state the result only outside these critical values.) 
Indeed, for $r\geq 2$ we have 
\[
C(\b_c(r)\mp 0,r) =(r-1)^{\frac{2(r-1)}{r(r-2)}}\times \left\{ 
\begin{array}{ll}
 r& 
\cr
 r (r-1)^\frac{r-2}{r}
& 
\end{array} \right.
\]
which jumps for $r\geq 3$. (For $r=2$ this expression has to be interpreted as the limit 
of the right hand side with $r\downarrow 2$.)



\medskip\noindent

The reader should notice the following: 
First of all we have shown the pointwise existence of the limit 
\[
(n_l)_{1\leq l\leq s} 
\mapsto \lim_{N\uparrow \infty}Q^N_{q, \beta, (r_1, \ldots, r_s)}
\bigl(k\bigl | (n_l)_{1\leq l\leq s} \bigr) \, .
\]
The notion of 
``continuity of limiting conditional probabilities''
that was introduced in Theorem \ref{thm:main_result_on_complete_graph} 
has the precise meaning of continuity 
of the right hand side 
as a function on the closed set $\PP$ of $s$-dimensional 
probability vectors with respect to the ordinary 
Euclidean topology. 
From the explicit limiting formula given in the theorem and the well-known 
knowledge of the jumps of the order parameter 
the proof of the first three parts of our main theorem 
 \ref{thm:main_result_on_complete_graph}  is now immediate.

 





\medskip\noindent
{\bf Proof of Theorem 
{\ref{thm:main_result_on_complete_graph} (i),(ii),(iii):}}
The points of discontinuity are precisely given by the 
values $n_k=\frac{\b_c(r_k)}{\b}$ for those $k$ with $r_k\geq 3$
for which $\frac{\b_c(r_k)}{\b}<1$. So (i) is immediate. 
To see (ii) and (iii) we use that $\b_c(r)$ 
is an increasing function of $r$. 
$\Cox$

\subsection{Typicality of continuity points -- ``almost sure Gibbsianness''}

What can be said about the measure of the discontinuity points? 
We will answer this question now and prove the remaining 
part (iv) of Theorem  \ref{thm:main_result_on_complete_graph}.  
To start with, from Theorem \ref{thm:Empirical measures} follows trivially 
by ``contraction'' that the  typical values of the order parameter in the fuzzy model are as follows. (Recall that $e_{l}$ 
is the unit vector in the $l$'th coordinate direction of $\R^{s}$.)
\begin{cor}  \label{cor:2.15}
We have 
\begin{eqnarray*}
\lefteqn{\lim_{N\uparrow\infty}
\mu^N_{q, \beta, (r_1, \ldots, r_s)}\Bigl(
\frac{1}{N}\sum_{x=1}^N (I_{\{Y(x)=1\}},\dots, I_{\{Y(x)=s\}})
\in \cdot \Bigr)} \\ 
&=& 
\left\{ 
\begin{array}{ll}
\d_{\frac{1}{q}(r_{1},r_{2},\dots,r_{s})}& \,\hbox{if}\,\,  \b<\b_c(q)
\cr
\l_0(q)\d_{\frac{1}{q}(r_{1},r_{2},\dots,r_{s})}+
\sum_{l=1}^{s} \frac{(1-\l_0(q)) r_{l}}{q}  
\d_{u(\b,q) e_{l}+\frac{1-u(\b,q)}{q}(r_{1},r_{2},\dots,r_{s})}
& \,\hbox{if} \,\, \b=\b_c(q)\cr
 \sum_{l=1}^{s} \frac{r_{l}}{q}  
\d_{u(\b,q) e_{l}+\frac{1-u(\b,q)}{q}(r_{1},r_{2},\dots,r_{s})}
& \,\hbox{if} \,\, \b>\b_c(q) \, .
\end{array}\right.
\end{eqnarray*}
\end{cor}
In other words, the values for the fuzzy densities $n_{l}$ that 
occur with non-zero probability are: The values $r_l /q$ 
in the high-temperature regime (including the critical point)
and the two values 
\[
n^{+}(\b,q,r_{l}) \equiv u(q,\b)+\frac{1-u(q,\b)}{q} r_{l}
\]
and
\[
n^{-}(\b,q,r_{l}) \equiv \frac{1-u(q,\b)}{q} r_{l}
\quad \Bigl( \leq n^{+}_{l}(\b,q,r_{l})\Bigl)
\]
in the low temperature regime (including the critical point). 

Now, the non-trivial question is: 
Can it happen that these values coincide with the points of discontinuity
of the limiting conditional probability, for certain choices of the parameter? 

The following proposition tells us that this can never be the case, 
and so the points of discontinuity are always atypical. This immediately 
proves (iv) of Theorem \ref{thm:main_result_on_complete_graph}.
As we will see the proof of the proposition is elementary 
but slightly tricky; it makes use 
of specific properties of the solution of the mean-field equation. 
In that sense it is the most difficult part of our 
analysis of the mean field fuzzy Potts model. 
\begin{prop}  \label{prop:300} Assume that $q>r\geq 2$. 
\begin{description}
\item{\bf (i)} For the high-temperature range $\b\leq  \b_{c}(q)$ we have 
\[
\frac{r}{q} < \frac{\b_{c}(r)}{\b} \, .
\]
\item{\bf (ii)} For the low-temperature range $\b\geq \b_{c}(q)$ we have that 
\[
n^{-}(\b,q,r)< \frac{\b_{c}(r)}{\b} <n^{+}(\b,q,r) \, .
\]
\end{description}
\end{prop}

\noindent{\bf Remark:}
(i) says that that the typical density of each fuzzy class is too small 
to create a first order transition. 
The left inequality of 
(ii) says that the typical density of a fuzzy class not 
containing the predominant  
spin-value of the underlying Potts model 
is  always too small to create a first order transition. The corresponding 
conditional Potts model is always in a high-temperature state. 
The right inequality of 
(ii) says that the typical density of the fuzzy class that 
contains the predominant 
spin-value of the underlying Potts model is always too big to 
create a first order transition. 
The corresponding 
conditional Potts model is always in a low-temperature state.  

\medskip\noindent
{\bf Proof:} 
The claim (i) follows from that fact that 
$\frac{r}{q}< \frac{\b_{c}(r)}{\b_c(q)}$ for all $q>r$. 
This in turn is implied by the fact that $\frac{\b_{c}(q)}{q}$ is decreasing 
in $q$. It is obvious that this holds for large enough $q$, by the 
explicit expression for $\b_{c}(q)$. It is elementary to verify that 
it holds in fact for any $q\geq 2$. 

Next we prove (ii).
We show first the right inequality
which is equivalent to 
\[
u(q,\b)> \frac{q }{q-r}\frac{\b_{c}(r)}{\b}-\frac{r}{q-r} \, .
\]
By Theorem \ref{thm:Empirical measures}
the order parameter $u(q,\b)$ is an increasing function 
in $\b$.    The right hand side is decreasing in $\b$. 
So it suffices to prove the inequality for $\b=\b_c(q)$. 
Using $u(q,\b_c(q))=\frac{q-2}{q-1}$ 
this can be put equivalently as 
\begin{equation}\label{eq:3.4}
\b_{c}(r)< 
\b_c(q)\left(1-\frac{q-r}{q(q-1)}\right) \, .
\end{equation}
We will use now the elementary property that
\begin{eqnarray} \label{eq:3.5}
\b_c(q)<q, \quad \hbox{for all real }q> 2 \, .
\end{eqnarray}
This implies also that $\b_c(q)$ is concave because 
\[
\b''_c(q)=\frac{-2 q (q-2)+ 4 (q-1)\log(q-1)}{(q-2)^3 (q-1)}
\]
and the denominator is negative, by the last inequality. 

In order to show (\ref{eq:3.4}) we note, by concavity that 
\begin{eqnarray}\label{eq:3.7}
\b_{c}(r)\leq \b_c(q)+\b'_c(q) (r-q)
\end{eqnarray}
and show that the right hand side 
of (\ref{eq:3.7}) is strictly bounded from above by 
the right hand side of (\ref{eq:3.4}). 
But the latter statement is equivalent to 
\[
\b'_{c}(q)>
\b_c(q)\frac{1}{q(q-1)} \, .
\]
Computing the logarithmic derivative 
$\frac{\b'_{c}(q)}{\b_c(q)}$ we see that this is equivalent to 
\[
\frac{1}{q-1}-\frac{1}{q-2}+\frac{1}{(q-1)\log(q-1)}>
\frac{1}{q(q-1)} \, .
\]
This inequality in turn reduces after trivial computation 
to the statement (\ref{eq:3.5}) and this concludes the proof of 
the right inequality of (i).

Let us come to the proof of the left inequality of (ii). 
The claim says 
$\frac{1-u(q,\b)}{q} r< \frac{\b_{c}(r)}{\b}$. 
Using the mean-field equation we may write 
\[
1-u(q,\b)
=\frac{q}{e^{+\b u(q,\b)}+q-1} \, .
\]
So the claim is equivalent to 
\[
\b \frac{r}{\b_c(r)}<e^{+\b u(q,\b)}+q-1 \, .
\]
Now, the left hand side 
is increasing as a function of $r$, for $r\geq 2$. So the claim 
follows from 
\[
\b \frac{q-1}{\b_c(q-1)}<e^{+\b u(q,\b)}+q-1 
\]
for all $q\geq 3$.
Next we use again 
that the order parameter $u(q,\b)$ is an increasing function of $\b$. 
Thus the last inequality follows if we can show 
\[
\b_c(q) \frac{q-1}{\b_c(q-1)}<e^{+\b_c(q) u(q,\b_c(q))}+q-1 \, .
\]
We have $e^{+\b_c(q) u(q,\b_c(q)+0)}=(q-1)^2$ from the explicit expressions and
so the last inequality is equivalent to 
\[
\frac{\b_c(q)}{\b_c(q-1)}<q \, . 
\]
It is elementary to verify from the explicit expression for $\b_c(q)$ that 
this actually holds for all $q\geq 3$. 
$\Cox$

\end{document}